\documentclass[12pt]{amsart}
\usepackage{amssymb,latexsym}
\usepackage{amsrefs}

\swapnumbers

\theoremstyle{plain}
\newtheorem{blank}{}[section]

\newtheorem*{thma}{Theorem~A}
\newtheorem*{thma0}{Theorem~A${_0}$}
\newtheorem*{thmb}{Theorem~B}
\newtheorem*{thmc}{Theorem~C}

\newtheorem*{cor}{Corollary}

\theoremstyle{definition}
\newtheorem{dblank}[blank]{}

\DeclareMathOperator{\Cl}{cl}
\DeclareMathOperator{\dis}{d}
\DeclareMathOperator{\fsig}{{F}_{\sigma}}
\DeclareMathOperator{\Fr}{fr}

\DeclareMathOperator{\graph}{graph}
\DeclareMathOperator{\Int}{int}
\newcommand{\cf}{\textit{cf.}}

\newcommand{\tx}{\textup}
\newcommand{\lqu}{\textup{``}}
\newcommand{\rqu}{\textup{''}}

\newcommand{\abs}[1]{\left\lvert#1\right\rvert}
\newcommand{\set}[1]{\{\, #1\,\}}
\newcommand{\Set}[1]{\left\{\, #1\,\right\}}
\newcommand{\olim}[1]{\lim_{#1\to 0^+}}
\newcommand{\ilim}[1]{\lim_{#1\to +\infty}}
\newcommand{\inv}{^{-1}}
\newcommand{\rest}{{\restriction}}
\newcommand{\e}{\epsilon}
\newcommand{\N}{\mathbb N}
\newcommand{\Q}{\mathbb Q}
\newcommand{\R}{\mathbb R}
\newcommand{\Z}{\mathbb Z}
\newcommand{\rbar}{\overline\R}
\newcommand{\rr}{\mathfrak R}
\newcommand{\nbd}{\nobreakdash-\hspace{0pt}}
\newcommand{\omin}{o\nbd minimal}
\newcommand{\dfbl}{definable}
\newcommand{\odfbl}{$\emptyset$\nbd\dfbl}

\newcommand{\bsig}{\boldsymbol\Sigma}
\newcommand{\bpi}{\boldsymbol\Pi}
\newcommand{\nul}{$\rr$\nbd null}

\makeatletter
   \def\th@plain{\slshape}
   \makeatother

\begin{document}

\title[Expansions of the real field by open sets]
{Expansions of the real field by open sets:\\
definability versus interpretability }

\author[H. Friedman]{Harvey Friedman}
\address
{Department of Mathematics\\
The Ohio State University\\
231 West~18th Avenue\\
Columbus, Ohio 43210, USA}
\email {friedman@math.ohio-state.edu}
%\urladdr{www.math.ohio-state.edu/\textasciitilde friedman}

\author[K. Kurdyka]{Krzysztof Kurdyka}
\address
{Laboratoire de Math\'ematiques (LAMA)\\ Universit\'e de
  Savoie\\UMR 5127 CNRS\\  73376 Le Bourget-du-Lac cedex
  France}
\email {Krzysztof.Kurdyka@univ-savoie.fr}
%\urladdr{http://www.lama.univ-savoie.fr/sitelama/Membres/pages_web/KURDYKA}

% Corresponding author
\author[C. Miller]{Chris Miller}
\address
{Department of Mathematics\\
The Ohio State University\\
231 West~18th Avenue\\
Columbus, Ohio 43210, USA}
\email{miller@math.ohio-state.edu}
%\urladdr{http://www.math.ohio-state.edu/\textasciitilde miller}

\author[P. Speissegger]{Patrick Speissegger}
\address
{Department of Mathematics \& Statistics \\
  McMaster University 1280 Main Street West Hamilton, Ontario L8S 4K1,
  Canada}
  \email{speisseg@math.mcmaster.ca}
%\urladdr{http://www.math.mcmaster.ca/\textasciitilde speisseg}

%\thanks
%{The first author was supported by NSF Grant No.~DMS-734164.

\begin{abstract}
An open $U\subseteq \mathbb R$ is produced such that $(\mathbb
R,+,\cdot,U)$ defines a Borel isomorph of $(\mathbb
R,+,\cdot,\mathbb N)$ but does not define $\mathbb N$. It follows
that $(\mathbb R,+,\cdot,U)$ defines sets in every level of the
projective hierarchy but does not define all projective sets. This
result is elaborated in various ways that involve geometric measure
theory and working over o-minimal expansions of $(\mathbb
R,+,\cdot\,)$. In particular, there is a Cantor set $K\subseteq
\mathbb R$ such that for every exponentially bounded o-minimal
expansion $\mathfrak R$ of $(\mathbb R,+,\cdot\,)$, every subset of
$\mathbb R$ definable in $(\mathfrak R,K)$ either has interior or is
Hausdorff null.
\end{abstract}

\thanks{\today. Some version of
this document has been submitted for publication. Comments are
welcome. Miller is the corresponding author.}



%\thanks{*Corresponding author}

%\thanks{\today}

\maketitle

The reader is assumed to be familiar with the basics of first-order
definability over the real field $\rbar:=(\R,+,\cdot\,)$, especially
\omin ity. Requisite material can be found in van den Dries and
Miller~\cite{geocat}. We refer to Kechris~\cite{kechris} and
Mattila~\cite{mattila} for basic descriptive set theory and
geometric measure theory. We say that a subset of $\R^n$ is
\textbf{constructible} if it is a boolean combination of open
subsets of $\R^n$. By Dougherty and Miller~\cite{DoM}, every
constructible $E\subseteq \R^n$ is a boolean combination of open
sets that are definable in $(\R,<,E)$; we use this fact without
further mention.

We begin with a simply-stated question: What can be said about the
sets that are definable, allowing arbitrary real parameters, in an
expansion $\rr$ of $\rbar$ by a collection of constructible subsets
of $\R$? First, every quantifier-free definable set is
constructible, hence Borel. Next, every existentially definable set
is $\bsig_1^1$ (also known as Souslin, Suslin, or analytic), and
every universally definable set is $\bpi_1^1$ (also known as
co-analytic). Continuing in this fashion, every definable set is
projective in the sense of descriptive set theory. All real
projective sets are definable in $(\rbar,\N)$~\cite{kechris}*{37.6}.
As a result, there are now many examples known where $\rr$ defines
all real projective sets; see~\cites{tameness, aph} for some
non-obvious ones. On the other hand, there are now many examples
known where every definable set is Borel; see~\cites{ominsparse,
fast,tameness,itseq} for some non-\omin\ ones. Heretofore, no other
behaviors have been documented. We show in this paper that there is
at least one other possibility: $\rr$ can define sets in every
$\bsig_{n+1}^1\setminus \bsig_n^1$, and thus of every projective
level, yet not define all real projective sets. As Borel
isomorphisms preserve projective level, this is immediate from

\begin{thma}
There is a closed $E\subseteq\R$ such that\tx:
\begin{enumerate}
\item[(a)]
$(\rbar,E)$ defines a Borel isomorph of $(\rbar,\N)$.
\item[(b)]
Every unary \tx(that is, contained in $\R$\tx) set definable in
$(\rbar,E)$ either has interior or is nowhere dense.
\end{enumerate}
\end{thma}

Expansions of $\rbar$ in which every unary definable set either has
interior or is nowhere dense have a number of good
properties~\cite{tameness}, but we shall not dwell on this here. On
the other hand, this condition is not strong enough to rule out a
significant difference between interpretability and definability of
interesting algebraic objects.

\begin{cor}
With $E$ as in Theorem~A\tx, $(\rbar,E)$ defines no proper
nontrivial subgroups of $(\R,+)$\tx, nor any proper noncyclic
subgroups of $(\R^{>0},\cdot\,)$\tx, yet $(\rbar,E)$ interprets
every real projective set\tx, in particular\tx, every projective
subfield of~$\rbar$.
\end{cor}
\noindent As we shall see (Theorem~B, below), we can choose $E$ so
that part~(b) of Theorem~A holds even for $(\rbar,\exp,E)$, in which
case neither does $(\rbar,E)$ define any proper nontrivial subgroups
of $(\R^{>0},\cdot\,)$.

We postpone beginning the proof proper of Theorem~A, but we outline
some of the main ideas now. In order to satisfy part~(b), it
suffices by~\cite{ominsparse}*{Theorem~A} and cell decomposition to
produce a closed $E\subseteq \R$ such that:
\begin{itemize}
\item[(b)\textprime]
For every $n\in\N$, bounded open semialgebraic cell $U\subseteq
\R^n$, and bounded continuous semialgebraic $f\colon U\to\R$, the
image $f(U\cap E^n)$ is nowhere dense, where $E^n$ is the $n$-th
cartesian power of~$E$.
\end{itemize}
Given any uncountable $\bsig_1^1$ set $E\subseteq \R$ such
that~(b)\textprime\ holds, there is a Cantor set\footnote{For our
purposes, a \textbf{Cantor set} is a subset of~$\R$ that is
nonempty, compact, and has neither interior nor isolated points.}
$K$ such that both~(a) and~(b)\textprime\ hold with $E$ replaced by
$K$; this is fairly easy modulo known descriptive set theory and
some definability tricks (\ref{newcantor}, below). Thus, it suffices
to find a Cantor set $E$ such that condition~(b)\textprime\ holds
for $E$. In order to motivate further developments, we consider a
naive approach that we could not make work. Let $E$ be a Cantor set
such that every $E^n$ is Hausdorff null (that is, has Hausdorff
dimension zero). By cell decomposition, we reduce further to the
case that $f$ is $C^1$ and nowhere locally constant. Write $U$ as
the union of the compact sets~$A_r$, $r>0$, where $A_r$ is the set
of $x\in U$ whose distance to the boundary of $U$ is at least $1/r$.
Each restriction $f\rest A_r$ is continuous and Lipshitz, so
$f(E^n\cap A_r)$ is compact and Hausdorff null. Since $f$ is nowhere
locally constant, $f(E^n\cap A_r)$ is also Cantor for all
sufficiently large $r$. Hence, $f(U\cap E^n)$ is the union of a
``semialgebraically parameterized'' increasing family of Hausdorff
null Cantor sets. But we see no way to conclude from this that
$f(U\cap E^n)$ is nowhere dense, which is required in order to
employ the aforementioned technology from~\cite{ominsparse}.  The
fundamental shortcoming of this approach is that it appears not to
account for how limit points of $f (U\cap E^n)$ are formed at the
boundary of $U$. We overcome this by a more careful choice of $E$
based on an analysis of the behavior of bounded semialgebraic
functions near their points of discontinuity. It turns out to be
just as easy to work in the more general setting of \omin\
expansions of $\rbar$ and prove stronger statements. In doing so, we
establish some results in \omin ity (\ref{bddcell}
through~\ref{Scellsexp}) that seem to be new to the literature even
as semialgebraic or subanalytic geometry. We also prove a
generalization and some variants of Theorem~A, one of which we state
now, leaving others for later.

A structure on $\R$ is \textbf{exponentially bounded} if for every
definable $f\colon \R\to\R$ there exists $m\in\N$ such that $f$ is
bounded at $+\infty$ by the $m$\nbd th compositional iterate
$\exp_m$ of $\exp$. It is an easy consequence of quantifier
elimination that $\rbar$ itself is \omin\ and exponentially bounded.
%(indeed, polynomially bounded).

\begin{thmb}
There is a Cantor set $K$ such that $(\rbar,K)$ defines a Borel
isomorph of $(\rbar,\N)$ and\tx, for every exponentially bounded
\omin\ expansion $\rr$ of $\rbar$\tx, every unary set definable in
$(\rr,K)$ either has interior or is Hausdorff null.
\end{thmb}

For the following reasons, we regard Theorem~B as a natural
extension of Theorem~A and its original motivating question. (i)~By
cell decomposition, every \omin\ expansion of the real line $(\R,<)$
is interdefinable with a structure on $\R$ generated by a collection
of open $U_{\alpha}\subseteq\R^{n(\alpha)}$, $\alpha$ ranging over
some index set. (ii)~For any expansion of $(\R,<)$, if every unary
definable set either has interior or is Hausdorff null, then every
unary definable set either has interior or is nowhere dense. (For
every $A\subseteq \R$ and open interval $I$, at most one of $I \cap
A$ and $I\setminus A$ is Hausdorff null.) (iii)~By growth
dichotomy~\cite{exphard}, Pfaffian closure~\cite{pfaffclosure}, and
Lion \textit{et al.}~\cite{lms}, if $\rr$ is an exponentially
bounded \omin\ expansion of $\rbar$, then so is $(\rr,\exp)$. In
particular, $(\rbar,\exp)$ is \omin\ and exponentially bounded.
(iv)~There are now many examples of expansions of $\rbar$ that are
known to be exponentially bounded and \omin.\footnote{There are also
exponentially bounded expansions of $\rbar$ that are not \omin. Such
structures are closely related to \omin ity\ in a certain way, but
cannot be obtained by expanding $\rbar$ by collections of
constructible sets (of any arities), and thus are irrelevant for
present purposes. See~\cite{opencore} for details.}

It is easy to see that a Cantor set is interdefinable over $\rbar$
with the set of midpoints of its complementary intervals
(see~\ref{discretecor}), so Theorem~B holds with ``discrete'', even
``countable'', in place of ``Cantor'', with similar modifications to
Theorem~A and its corollary.

As of this writing, every expansion of $\rbar$ known to be \omin\ is
exponentially bounded, thus partly justifying our decision to
postpone the statement of our most general version (Theorem~C in
\S\ref{S:thmb}). Another reason is simply to avoid for now having to
introduce further notation and technical definitions.

Though one might be tempted to regard both Theorems~A and~B as
teratological, we believe that the techniques of the proofs are
interesting in their own right, and potentially useful for other
settings.

Here is an outline of the remainder of this paper. We begin in
Section~\ref{S:prelim} with some preliminaries, including some
results in pure \omin ity that we believe will be useful in other
settings. We prove Theorem~B in Section~\ref{S:thmb}, as well as
some variants and corollaries. We close in Section~\ref{S:conc} with
discussion and open issues.

\section{Preliminaries}\label{S:prelim}

We begin by establishing some global conventions and notation.
Throughout, ``definable'' (in some first-order structure) means
``definable with parameters'', while ``\odfbl'' means ``definable
without parameters''. The variables $j,k,m,n$ range over $\N$, the
non-negative integers. Given a set $A$, its $n$\nbd fold cartesian
power is denoted by $A^n$, with $A^0:=\{0\}$. Whenever convenient,
we identify $A^m\times A^n$ with $A^{m+n}$, in particular,
$A^m\times A^0\cong A^0\times A^m\cong A^{m}$.
%Given $C\subseteq A^{m+n}$ and $x\in A^m$, we let $C_x$
%denote the \textbf{fiber} of $C$ over $x$, that is, $C_x=\set{y\in
%A^n: (x,y)\in C}$.
If $A$ belongs to a topological space, we denote its interior by
$\Int(A)$, closure by $\Cl(A)$,
%boundary by $\Bd(A)$ ($=\Cl(A)\setminus \Int(A)$),
and frontier by $\Fr(A):=\Cl(A)\setminus A$. If $A\subseteq \R^n$,
then all of these sets (in the usual topology) are \odfbl\ in
$(\R,<,A)$. Given a set $B$, we identify a function $f\colon A^0\to
B$ with the constant $f(0)\in B$. Given a function $f\colon A\to B$
and $A'\subseteq A$, we let $f\rest A'$ denote the restriction of
$f$ to $A'$. Limits of functions are always taken with respect to
the declared domain of the function, and similarly with limits
superior and inferior. Metric notions are taken with respect to the
sup norm $\abs{x}:=\sup\{\abs{x_1},\dots,\abs{x_n}\}$. In
particular, for $x\in \R^n$ and $r>0$, put $B(x,r)=\set{y\in\R^n:
\abs{x-y}<r}$.

Next are some crucial technical results in \omin ity that seem to be
new to the literature even as semialgebraic or subanalytic
geometry.\footnote{We would appreciate any information to the
contrary.}

\begin{blank}\label{bddcell}
Let $\rr$ be an \omin\ expansion of $(\R,<)$\tx, $C\subseteq \R^n$
be a bounded open cell\tx, and $f\colon C\to \R$ be definable\tx,
continuous\tx, and bounded. Then there is a definable $X\subseteq
\Fr(C)$ such that $\dim(\Fr(C)\setminus X)<n-1$ and $f$ extends
continuously to $C\cup X$.
\end{blank}

(Recall that $\dim\emptyset=-\infty$ by convention.)

\subsubsection*{Note}
The resulting extension of $f$ to $C\cup X$ is necessarily
definable, as it is given by $y\mapsto \lim_{x\to y}f(x)$.

\begin{proof}
As is often the case in o-minimality, we find it convenient to prove
simultaneously a related condition. We proceed by induction on
$n\geq 1$ to show the following in turn.
\begin{itemize}
\item[(i$_n$)]
There is a definable $Y\subseteq \Fr(C)$ such that
$\dim(\Fr(C)\setminus Y)<n-1$ and $C$ is locally connected at every
$y\in Y$.
\item[(ii$_n$)]
There is a definable $X\subseteq \Fr(C)$ such that
$\dim(\Fr(C)\setminus X)<n-1$ and $f$ extends continuously to $C\cup
X$.
\end{itemize}
(A subset $A$ of a topological space $X$ is \textbf{locally
connected} at $x\in X$ if for every open neighborhood $U$ of $x$
there is an open neighborhood $V\subseteq U$ of $x$ such that $V\cap
A$ is connected.)

If $n=1$, then $C$ is an open interval, so the result is immediate
from the monotonicity theorem. Let $n\geq 1$ and assume the result
for~$n$. Let $C\subseteq \R^{n+1}$ be a bounded open cell.

(i$_{n+1}$). Let $D$ be the projection of $C$ on the first $n$
variables. Then $D$ is a bounded open cell and there exist bounded
definable continuous functions $g,h\colon D\to \R$ such that $g<h$
and $C=\set{(x,r): x\in D\And g(x)<r<h(x)}$. Inductively, there
exist definable $Z\subseteq \Fr(D)$ and definable continuous
$G,H\colon D\cup Z\to\R$ such that $\dim(\Fr(D)\setminus Z)<n-1$,
$D$ is locally connected at every $z\in Z$, $g=G\rest D$, and
$h=H\rest D$. Put
$$
Y=\graph(g)\cup\graph(h)\cup\set{(z,r): z\in Z\And G(z)<r<H(z)}.
$$
Then
$$
\Fr(C)\setminus Y\subseteq \Fr(\graph(g))\cup\Fr(\graph(h))\cup
[(\Fr(D)\setminus Z)\times \R],
$$
so $ \dim(\Fr(C)\setminus Y)<n$. We now show that $C$ is locally
connected at every $y\in Y$. As $D$ is open and cells are connected,
$C$ is locally connected at every point of $\graph(g)\cup\graph(h)$.
Let $z\in Z$ and $r\in \R$ be such that $G(z)<r<H(z)$. Let $U$ be an
open set containing $(z,r)$. We must find an open box about the
point $(z,r)$ that is contained in $U$ and whose intersection with
$C$ is connected. By continuity of $G$ and $H$, there is an open box
$B\times I\subseteq U$ about $(z,r)$ such that $(B\times I)\cap C$
is disjoint from $\graph(g)\cup\graph(h)$. Since $D$ is locally
connected at $z$, we may shrink $B$ so that $B\cap D$ is connected.
Then $(B\times I)\cap C$ is connected, because $(B\times I)\cap
C=(B\cap D)\times I$.

(ii$_{n+1}$). Let $f\colon C\to \R$ be definable, continuous, and
bounded. Let $Z$ be the set of all $z\in \Fr(C)$ such that
$\lim_{x\to z}f(x)$ exists. We claim that $\dim(\Fr(C)\setminus
Z)<n$. Suppose not. Then there is a cell $E\subseteq \Fr(C)$ such
that $\dim E=n$ and $\liminf_{x\to y}f(x)<\limsup_{x\to y}f(x)$ for
every $y\in E$. By (i$_{n+1}$), we may shrink $E$ so that $C$ is
locally connected at every $y\in E$. Then
$$
\Fr(\graph(f))\supseteq \set{(y,r):y\in E\And \liminf_{x\to
y}f(x)<r<\limsup_{x\to y}f(x)},
$$
yielding the absurdity
%\begin{align*}
%\dim \Fr(\graph(f))&\geq \dim E +1\\
%&=n+1\\
%&=\dim \graph(f)\\
%&>\dim \Fr(\graph(f)).%,\qquad\text{(\cite{ominbook}*{pg.~67}.)}
%\end{align*}
$$
\dim \Fr(\graph(f))\geq \dim E +1=n+1=\dim \graph(f)>\dim \Fr(\graph(f)).
$$
(See van den Dries~\cite{ominbook}*{pp.~65,~68}.) Define $g\colon
Z\to\R$ by $g(z)=\lim_{x\to z}f(x)$. Let $X$ be the set of points of
continuity of $g$. By cell decomposition, $\dim(Z\setminus X)<\dim
Z$, so $\dim(\Fr(C)\setminus X)<n$. Finally, define $h\colon C\cup
X\to\R$ by $h\rest C=f$ and $h\rest X=g\rest X$. Then $h$ is
continuous and extends $f$.
\end{proof}

\subsubsection*{Remarks}
%(i)~The result holds with ``\odfbl'' instead of ``definable''.
(i)~The result and its proof hold for all abstract \omin\ structures
(as defined in~\cite{ominbook}) provided that the definition of
locally connected is relativized to definable connectedness. (ii)~If
$C=\set{(x,y)\in\R^2: 0<y<x<1} $ and $f(x,y)=y/x$, then $f$ is
bounded and does not extend continuously to the origin. (iii)~If $
C=\set{(x,y,z)\in\R^3: \abs{x}<1\And 0<y<1\And
-1<z<\sqrt{\abs{x}}/y}, $ then $C$ is bounded and not locally
connected at any point of $\set{(0,0,z): 0<z<1}$.

We need yet more refined results, for which we require some
technical
definitions and notation. %In order to motivate them, let us consider
%in more detail the example in Note~(ii) above. While $f$ does not
%extend continuously to the origin, for every set of the form
%$S:=\set{(x,y): 0<y<x^\alpha\leq b}$ where $0<b<1<\alpha$, not only
%does $f\rest S$ extend continuously to the origin, but to all of
%$\Cl(S)$. Moreover, the same is true of the set $\set{(x,y):
%0<x-x^\alpha<y<x\leq b}$.

%\begin{dblank}
Define the \textbf{corners} in $\R^n$ inductively as
follows.\footnote{More precisely, we define the corners ``at
$0^+$''. Corners can be defined relative to any point in
$[-\infty,\infty]^n$ regarded with a fixed appropriate sign
condition, but we shall not need any of these variants in this
paper.} (i)~$\R^0$ is the only corner in $\R^0$. (ii)~If $C\subseteq
\R^n$ is a corner and $f\colon C\to (0,\infty)$ is continuous, then
$$\set{(x,t)\in\R^{n+1}: x\in C\And 0<t<f(x)}$$ is a corner in
$\R^{n+1}$. We note some easy facts.
%\begin{enumerate}
%\item
Every corner in $\R^n$ is an open cell contained in $(0,\infty)^n$.
%\item
%For every corner $C$ in $\R^n$, there is a corner $C'$ in $\R^n$
%such that $\Cl(C')\cap (0,\infty)^n\subseteq C$.
%\item
The projection on the first $m$ coordinates of a corner in
$\R^{m+n}$ is a corner in $\R^m$.
%\item
%For every corner $C$ in $\R^{n+1}$ there exists $\e>0$ such that
%$(0,\e)\times \{0\}^{n}\subseteq \Cl(C)$.
%\item
%For every corner $C$ in $\R^n$ there is a corner $C'$ in $\R^n$ such
%that $\Cl(C')\cap (0,\infty)^n\subseteq C$.
%\item
For every cell decomposition $\mathcal C$ of $\R^n$ that is
compatible with $(0,\infty)^n$ there is a unique $C\in\mathcal C$
such that $C$ is a corner.
%\end{dblank}

\begin{dblank}\label{phidef}
We now define some special corners. Let $\Phi$ be the collection of
all homeomorphisms of $[0,\infty)$. Let $\mathcal S_n$ be the
collection of all nonempty sets of the form
\begin{equation*}%\label{Sdef}
\set{x\in (0,\infty)^n:
\phi_n(x_n)<\phi_{n-1}(x_{n-1})<\dots<\phi_{1}(x_{1})<b}
\end{equation*}
where $b\in (0,+\infty]$ and $\phi_1,\dots,\phi_n\in\Phi$. Note that
this set is an open cell of the structure
$(\R,<,\phi_1,\dots,\phi_n)$. We interpret $\mathcal S_0$ as $\R^0$.
Let $S\in\mathcal S_n$. Note the easy facts: $S$ is a corner; for
every $m\leq n$, the projection of $S$ on the first $m$ coordinates
belongs to $\mathcal S_m$; and $\Cl(S)\setminus
[(0,\infty)^{n-1}\times \R]=\{0\}^n$. With $b=+\infty$ and
$\phi_i=t\mapsto 2^{i}t$ for $i=1,\dots,n$, we obtain the set
%$$
%\mathbb S_n:=\set{x\in \R^n: 0<2^{n-1}x_n<\dots<2x_2<x_{1}}.
%$$
$$
\mathbb S_n:=\Set{x\in \R^n:
0<x_n<\frac{x_{n-1}}{2}<\dots<\frac{x_{1}}{2^{n-1}}}.
$$
Note that $\mathbb S_n$ is definable in $(\R,<,+)$.
\end{dblank}

\begin{blank}\label{Scells}
Let $C\subseteq \R^n$ be a corner such that $(\R,<,+,C)$ is \omin.
Then there exists $S\in\mathcal S_n$ definable in $(\R,<,+,C)$ such
that $\Cl(S)\cap (0,\infty)^n\subseteq C$.
\end{blank}

\begin{proof}
We proceed by induction on $n\geq 1$. (The case $n=0$ is trivial.)
If $n=1$, then $C$ is an open interval $(0,s)$ for some $s\in
(0,\infty]$, so let $S$ be any $(0,r)$ with $r\in (0,s)$.

Let $n\geq 1$ and assume the result for $n$. Let $C\subseteq
(0,\infty)^{n+1}$ be a corner such that $(\R,<,+,C)$ is \omin. It
suffices to consider the case that $C$ is bounded. Let $D$ be the
projection of $C$ on the first $n$ coordinates. Then $D$ is a
bounded definable corner and there is a continuous definable
$f\colon D\to (0,\infty)$ such that
$$
C=\set{(x,x_{n+1})\in \R^{n+1}: x\in D\And 0<x_{n+1}<f(x)}.
$$
Inductively, there is a definable $S'\in\mathcal S_n$ such that $
\Cl(S')\cap (0,\infty)^n\subseteq D$. The projection of $S'$ on the
last coordinate is an open interval $(0,a)$, and for every $t\in
(0,a)$, the set $\set{x\in \Cl(S'): x_n=t}$ is compact. As $f$ is
continuous on $\Cl(S')\cap (0,\infty)^n$, for all $t\in (0,a)$ we
have $0<\inf \set{f(x): x\in \Cl(S')\And x_n=t}<\infty$. Define
$g\colon (0,a)\to \R$ by
$$
g(t)=\min(t, \inf\set{f(x)/2: x\in \Cl(S')\And x_n=t}).
$$
Note that $g$ is definable. As $g>0$ and $\olim g(t)=0$, there
exists by the monotonicity theorem some $a'\in (0,a]$ such that
$g\rest (0,a')$ is continuous and strictly increasing. The set
$\set{x\in S':x_n<a}$ is a corner, so inductively, we may shrink
$S'$ so that $g$ is continuous and strictly increasing. Put
$
S=\set{(x,x_{n+1}): x\in S'\And x_{n+1}<g(x_n)}.
$
Now,
$$
\Cl(S)\cap (0,\infty)^{n+1}\subseteq \set{(x,x_{n+1}): x\in
\Cl(S')\cap (0,\infty)^n\And 0<x_{n+1}\leq f(x)/2}\subseteq C,
$$
so it suffices to show that $S\in\mathcal S_{n+1}$. Extend $g$ to
$\phi\in\Phi$ by setting $\phi(0)=0$ and $\phi(t)=t-a+g(a)$ for
$t>a$. Write
$$
S'=\set{x\in (0,\infty)^n: \phi_n(x_n)<\dots<\phi_{1}(x_{1})<b}
$$
as in~\ref{phidef}. Then $\phi_n\circ\phi\inv\in\Phi$ and
\[
S=\set{(x,x_{n+1})\in (0,\infty)^{n+1}:\phi_n\circ\phi\inv
(x_{n+1})<\phi_n(x_n)<\dots<\phi_{1}(x_{1})<b}\in\mathcal
S_{n+1}.\qedhere
\]
\end{proof}

\begin{blank}\label{Scellscont}
Let $S\in\mathcal S_n$\tx, $f\colon S\to\R$ be bounded and
continuous\tx, and $(\R,<,+,f)$ be \omin. Then there exists
$S'\in\mathcal S_n$ definable in $(\R,<,+,f)$ such that $S'\subseteq
S$ and $f\rest S'$ extends continuously to $\Cl(S')$.
\end{blank}

\begin{proof}
%Assume that $n>0$. Let $\pi$ denote projection on the first $n-1$
%coordinates. By~\ref{bddcell} and~\ref{Scells}, we reduce to the
%case that there is a definable continuous $g\colon S\cup (\pi
%S\times \{0\})\to\R$ such that $f=g\rest S$. By~\ref{Scells}, we may
%assume that $g$ extends continuously to $\Cl(\pi S)\cap
%(0,\infty)^{n-1}$. By replacing $f$ with $x\mapsto f-g\circ
%(\pi,0)$, we reduce to the case that $\olim{t}f(\pi x,t)=0$ for all
%$x\in S$. It follows that the set $\set{x\in S: \abs{f(x)}\leq
%x_{n}}$ contains a corner, so by~\ref{Scells}, we reduce to the case
%that $\abs{f(x)}\leq x_{n}$ for all $x\in S$. Again by~\ref{Scells},
%we reduce to the case that $f$ extends continuously to $\Cl(S)\cap
%(0,\infty)^n$ and $\abs{f(x)}\leq x_{n}$ for all $x\in S$. Then
%$\lim_{x\to y}f(x)=0$ for all $y\in \Cl(\pi S)\times \{0\}$, and we
%are done.
We proceed by induction on $n\geq 1$. The case $n=1$ is immediate
from the monotonicity theorem. Let $n\geq 1$ and assume the result
for $n$. Let $\pi$ denote projection on the first $n$ coordinates.
By~\ref{bddcell} and~\ref{Scells}, we may shrink $\pi S$ (hence also
$S$) so that there is a definable continuous $g\colon S\cup (\pi
S\times \{0\})\to\R$ with $f=g\rest S$. Inductively, we reduce to
the case that $g\rest (\pi S\times \{0\})$ extends continuously and
definably to $\Cl(\pi S)\times \{0\}$; we denote this extension just
by $g$. By continuity, the set $\set{x\in S: \abs{f(x)-g(\pi
x,0)}\leq x_{n+1}}$ contains a corner, so we reduce by~\ref{Scells}
to the case that $\abs{f(x)-g(\pi x,0)}\leq x_{n+1}$ for all $x\in
S$. Hence, $\lim_{x\to (y,0)}f(x)=g(y,0)$ for all $y\in \Cl(\pi S)$,
so $f$ extends continuously to $\Cl(\pi S)\times \{0\}$. Finally, we
shrink $S$ again by~\ref{Scells} so that $f$ is continuous on
$\Cl(S)\cap (0,\infty)^{n+1}$. Then $f$ extends continuously to
$\Cl(S)$.
\end{proof}

%\begin{blank}\label{Scellscont}
%Let $f\colon (0,\infty)^n\to\R$ be bounded near $0$ and $(\R,<,+,f)$
%be \omin. Then there there exists $S\in\mathcal S_n$ definable in
%$(\R,<,+,f)$ such that $f\rest S$ is continuous and extends
%continuously to $\Cl(S)$.
%\end{blank}

\begin{dblank}\label{Tdef}
Let $\mathcal T_n$ be the group of symmetries, regarded as linear
transformations $\R^n\to\R^n$, of the polyhedron inscribed in the
unit ball in $\R^n$ whose vertices are the intersections of the unit
sphere in $\R^n$ with the set $\set{tu: t>0\And u\in\{-1,0,1\}^n}$.
For example, $\mathcal T _2$ is the symmetry group of the octagon.
With $\mathbb S_n$ as in~\ref{phidef}, we have $
\R^n=\bigcup_{T\in\mathcal T_n}T(\Cl(\mathbb S_n))=
\bigcup_{T\in\mathcal T_n}\Cl(T(\mathbb S_n))$. As an immediate
consequence of cell decomposition, \ref{Scells},
and~\ref{Scellscont}, we obtain a key technical lemma:
\end{dblank}

\begin{blank}\label{Scellscontmain}
Let $A\subseteq \R^n$ and $f\colon A\to\R$ be such that $(\R,<,+,
f)$ is \omin. Let $f$ be bounded near $y\in\Cl(A)$. Then there
exists $S\in\mathcal S_n$ definable in $(\R,<,+, f)$ such that for
every $T\in \mathcal T_n$ and $m\leq n$\tx, the restriction of $f$
to $A\cap (y+T(\pi_m S \times\{0\}^{n-m}))$ is continuous and
extends continuously to the closure\tx, where $\pi_m$ denotes
projection on the first $m$ coordinates.
\end{blank}

\subsubsection*{Remark}
The result is easily modified to hold for \omin\ expansions of
arbitrary ordered groups. As one might imagine, the result can be
generalized significantly if $\rr$ also expands an ordered
field---recall the triangulation theorem---but it takes a bit of
effort to make this precise. As we shall not need any of these of
these generalizations in this paper, we leave details to the
interested reader.

\begin{dblank}\label{Sdef}
We now define some special elements of $\mathcal S_n$. First, define
$\psi\colon [0,\infty)\to\R$ by
$$\psi(t)=\begin{cases}0,& t=0\\
e^{-1/t},& 0<t<1\\
t-1+e^{-1},& t\geq 1.
\end{cases}
$$
Note that $\psi\in\Phi$. For $l\in\Z$, let $\psi_l$ be the $l$-th
compositional iterate of $\psi$. Every $\psi_l$ is definable in
$(\rbar,\exp)$. For $l\geq 0$, $\psi_l(t)=1/\exp_l(1/t)$ for all
sufficiently small $t>0$, and $\psi_{-l}(t)=1/\log_{l}(1/t)$ for all
sufficiently small $t>0$, where $\log_l$ denotes the ultimately
defined $l$-th compositional iterate of $\log$. If $\rr$ is an
exponentially bounded expansion of $\rbar$, then for each definable
$f\colon (0,b)\to (0,\infty)$ such that $\olim{t}f(t)=0$ there is a
least $l\in\Z$ such that $f(t)\geq \psi_{l}(t)$ as $t\to 0^+$; if
$\rr$ is also \omin, then $f(t)<\psi_{l-1}(t)$ as $t\to 0^+$. Define
sets $S_{n,l}\in\mathcal
S_n$ by %by taking each $\phi_i$ (as in~\ref{phidef}) to be
%$\psi_{-il}$, equivalently,
%
%inductively by
%\begin{align*}
%S_{0,l}&=\{0\}\\
%S_{1,l}&=(0,\infty)\\
%S_{n+1,l}&=\set{(x,r): x\in S_{n,l}\And 0<r<\psi_l(x_{n})},\quad
%n>1.
%\end{align*}
$$
S_{n,l}=\set{x\in \R^n:
0<x_n<\psi_l(x_{n-1})<\dots<\psi_{(n-1)l}(x_1)}.
$$
Every $S_{n,l}$ is an open cell of $(\rbar,\exp)$. An easy induction
shows that if $S\in\mathcal S_n$ and $(\rbar,S)$ is \omin\ and
exponentially bounded, then there exist $\delta,l>0$ such that
$B(0,\delta)\cap S_{n,l}\subseteq S$. Hence,
by~\ref{Scellscontmain},
\end{dblank}

\begin{blank}\label{Scellsexp}
Let $A\subseteq \R^n$ and $f\colon A\to\R$ be such that $(\rbar,f)$
is \omin\ and exponentially bounded. Let $f$ be bounded near
$y\in\Cl(A)$. Then there exist $\delta,j>0$ such that for every
$T\in \mathcal T_n$ and $m\leq n$\tx, the restriction of $f$ to
$B(y,\delta)\cap A\cap (y+T(S_{m,j}\times\{0\}^{n-m}))$ is
continuous and extends continuously to the closure.
\end{blank}

\subsubsection*{Remark}
The above holds with ``polynomially'' instead of ``exponentially''
by using the sets
$$ \Set{x\in \R^m: 0<x_n<
x_{n-1}^l<\dots<x_{1}^{l^{m-1}}}
$$ instead of the $S_{n,l}$.

Our results so far have not required working over $\R$ in that they
are easily generalized to abstract \omin\ structures. This now
begins to change.

%\begin{dblank}
Given $A\subseteq \R^n$ and $r>0$, let $N(A,r)$ be the infimum of
all $k$ such that $A$ is covered by $k$\nbd many closed cubes of
side length~$r$. The set $A$ is \textbf{Minkowski
null}\footnote{Also known under several other names and equivalent
formulations.} if $\olim{r}r^{\e}N(A,r)=0$ for all $\e>0$. While
Minkowski nullity is not generally preserved under countable unions
or $C^1$ images, it is much better behaved than Hausdorff nullity in
some other ways. We extend this notion relative to expansions $\rr$
of $\rbar$ by defining $A\subseteq \R^n$ to be \textbf{\nul} if
$\olim{r}f(r)N(A,r)=0$ for every definable $f\colon \R\to\R$ such
that $\olim{r}f(r)=0$. As $\rbar$ defines all rational power
functions, \nul\ implies Minkowski null. If $\rr$ is polynomially
bounded, then \nul\ is the same as Minkowski null. We have some
other easy basic facts that will be used often.
%\end{dblank}

\begin{blank}\label{nullfacts}
\begin{enumerate}
\item
\nul ity is preserved under taking subsets\tx, closure\tx, finite
unions\tx, finite cartesian products\tx, and images under Lipshitz
maps.
\item
\nul\ sets are nowhere dense and Hausdorff null.
\item
Countable unions of \nul\ sets are Baire meager and Hausdorff null.
\item\label{hnull}
If every bounded unary definable set either has interior or is
\nul\tx, then every unary definable set either has interior or is
Hausdorff null.
\item\label{ebd}
If $\rr$ is o\nbd minimal and exponentially bounded\tx, then
$A\subseteq \R^n$ is \nul\ if and only if $\displaystyle\lim_{r\to
0^+}r\exp_m(N(A,r))=0$ for all $m$.
\end{enumerate}
\end{blank}

\begin{proof}
We give a sketch for~(\ref{ebd}) but leave the rest as exercises.

Suppose that $A$ is \nul. Then $\lim_{r\to
0^+}N(A,r)/\log_m(\e/r)=0$ for every $m\in\N$ and $\e>0$. Hence,
$N(A,r)<\log_m(\e/r)$ for all sufficiently small $r>0$. Then
$r\exp_m(N(A,r))<\e$ for all sufficiently small $r>0$.

Conversely, suppose that $\lim_{r\to 0^+}r\exp_m(N(A,r))=0$ for all
$m$. We need only show that $ \lim_{r\to 0^+}N(A,r)/\log_m(1/r)=0 $
for all $m$. By assumption,
$$\lim_{r\to 0^+}r\exp_{m+1}(N(A,r))=0, $$
so $N(A,r)<\log_{m+1}(1/r)$ for all sufficiently small $r>0$. Now note
that
\[\olim{r}\log_{m+1}(1/r)/\log_m(1/r)=0.\qedhere
\]
\end{proof}

\begin{blank}\label{holderimage}
Let $\rr$ be an \omin\ expansion of $\rbar$. If $A\subseteq \R^m$ is
\nul\tx, $B\subseteq \R^m$ is compact\tx, and $f\colon B\to\R^n$ is
continuous and definable\tx, then $f(A\cap B)$ is \nul.
\end{blank}

\begin{proof}
We may take $A\subseteq B$. By generalized H\"older
continuity~\cite{geocat}*{C.15}, there is a definable $\phi\in\Phi$
such that $\abs{f(x)-f(y)}\leq \phi(\abs{x-y})$ for all $x,y\in B$.
Then $N(f(A),\phi(r))\leq N(A,r)$ for all $r>0$. Let $g\colon
(0,\infty)\to \R$ be definable such that $\olim{r}g(r)=0$. Then
\begin{equation*}
\olim{r}g(r)N(f(A),r)=\olim{r}g(\phi(r))N(f(A),\phi(r))
=\olim{r}(g\circ \phi)(r)N(A,r)=0.\qedhere
\end{equation*}
\end{proof}

We next recall a result from~\cite{ominsparse} and some minor
variants. Given a structure $\rr$ on $\R$ and $Y\subseteq \R$, let
$(\rr,Y)^{\#}$ denote the expansion $(\rr,(X))$ of $\rr$, where $X$
ranges over all subsets of all cartesian powers of $Y$.

\begin{blank}\label{sparsefact}
Let $A\subseteq \R$ and $\rr$ be an \omin\ expansion of
$(\R,<,+,1)$. Let $B\subseteq \R$ have no interior and be definable
in $(\rr,A)^{\#}$. Then there exists $f\colon \R^n\to\R$ definable
in $\rr$ such that $B\subseteq \Cl(f(A^n))$. If $B$ is bounded\tx,
then $f$ can be taken to be bounded by replacing it with
$\max(N,\abs{f})$ for some $N$ such that $B\subseteq [-N,N]$. This
all holds with \lqu\odfbl\rqu\ in place of \lqu\dfbl\rqu.
\end{blank}

Thus, as \nul ity is preserved under taking closure, we have

\begin{blank}\label{FMcase}
Let $\rr$ be an \omin\ expansion of $\rbar$ and $A\subseteq \R$ be
such that $f(A^n)$ is \nul\ for all bounded definable
$f\colon\R^n\to\R$. Then every bounded unary set definable in
$(\rr,A)^{\#}$ either has interior or is \nul.
\end{blank}

We leave the proof of the following amusing result as an
exercise.\footnote{We are not aware of this appearing elsewhere. We
would appreciate any information to the contrary.}

\begin{blank}\label{Tlemma}
For every $\Q$\nbd linearly independent $A\subseteq \R$\tx, the
function $x\mapsto \sum_{i=1}^n2^{i-1}x_i\colon \R^n\to\R$ is
injective on $A^n$.
\end{blank}
\noindent

Next is a combination of basic definability and descriptive set
theory.

\begin{blank}\label{newcantor}
Let $A\subseteq \R$ be uncountable and $\bsig_1^1$. Then there is a
Cantor set $K$ such that\tx:
\begin{enumerate}
\item
$(\R,<,+,1,K)$ $\emptyset$\nbd defines a Borel isomorph of
$(\rbar,\N)$.
\item
For every bounded $f\colon \R^n\to\R$ there is a finite set
$\mathcal H$ of bounded functions \odfbl\ in $(\R,<,+,1,f)$ such
that $f(K^n)\subseteq \bigcup_{h\in\mathcal H}h(A^{n(h)})$.
\end{enumerate}
\end{blank}

\begin{proof}
There exists $l\in\Z$ such that $A\cap [l,l+1]$ is uncountable. By
translation, we take $A\subseteq [1,2]$. Every uncountable analytic
set contains a Cantor set, so we take $A$ to be Cantor. Every Cantor
set contains a $\Q$\nbd linearly independent Cantor set
(use~\cite{kechris}*{19.1}), so we also take $A$ to be $\Q$\nbd
linearly independent. Let $D$ be the result of removing from $A$ its
maximum and minimum, and all left endpoints of the complementary
intervals of $A$. By a classical construction---for example, Gelbaum
and Olmstead~\cite{GO}*{I.8.14}---there is a strictly increasing
bijection $g\colon \R\to D$ whose compositional inverse is
continuous. Hence, $g$ is a Borel isomorphism and a closed map. Put
$$
X=\set{\bigl(g(x),g(y),g(x+y), g(xy),g(\dis_{\N}(x))\bigr):x,y\in
\R}
$$
where $\dis_\N\colon \R\to\R$ denotes the distance function to $\N$.
Define $T\colon \R^5\to\R$ by $T(x)=\sum_{i=1}^52^{i-1}x_i$. Put
$K=A\cup T(\Cl(X))$. As $X$ has no isolated points, the same is true
of $\Cl(X)$, hence also of $T(\Cl(X))$. Since $\Cl(X)\subseteq A^5$
and $T\rest A^5$ is injective (\ref{Tlemma}), $T(\Cl(X))$ is compact
and has no interior. Thus, $T(\Cl(X))$ is Cantor, hence so is $K$.

(1). We have $T(A^5)\subseteq (2,\infty)$, so $A=K\cap [1,2]$ and
$(\R,<,+,1,K)$ $\emptyset$\nbd defines $A$, hence also $D$. As $g$
is a closed map and $T\rest A^5$ is injective, we have $X=D^5\cap
T\inv(K\setminus A)$. Thus, $(\R,<+,1,K)$ $\emptyset$\nbd defines
$X$. Observe that $(\R,X)$ $\emptyset$\nbd defines the image under
$g$ of $(\rbar,\dis_{\N})$, hence also that of $(\rbar,\N)$.

(2). Straightforward, but tedious to write up in detail. We
illustrate the point via the case $n=2$ and leave the rest to the
reader. Since $K\subseteq A\cup T(A^5)$, we have
$$
K^2\subseteq A^2\ \cup\ A\times T(A^5)\ \cup\ T(A^5)\times A\ \cup\
(T(A^5))^2.$$ Let $f\colon \R^2\to \R$ be bounded. Put $\mathcal
H=\{h_1,\dots,h_4\}$, where
\begin{align*}
h_1&=f(x_1)\\
h_2&=f(x_1,T(x_2,\dots,x_6))\\
h_3&=f(T(x_1,\dots,x_5),x_6)\\
h_4&=f(T(x_1,\dots,x_5),T(x_6,\dots,x_{10})).
\end{align*}
Then $ f(K^2)\subseteq h_1(A^2)\cup h_2(A^6)\cup h_3(A^6)\cup
h_4(A^{10})$.
\end{proof}

\subsubsection*{Remark}By using distance functions, the set $X$
is easily modified to encode over $(\R,<,+,1)$ any expansion of
$(\R,<,+,1)$ by finitely many closed sets.

\section{Main results}\label{S:thmb}

In this section, we prove Theorem~B, as well as some variants and
corollaries.

\begin{proof}[Proof of Theorem~B]
By diagonalization, we fix a sequence $(r_k)$ of positive real
numbers such that $r_0=1$, $2r_{k+1}<r_{k}$ for all $k$, and
$\ilim{k}r_{k+1}/\psi_j(r_k)=0$ for all $j$ ($\psi_j$ as
in~\ref{Sdef}). Define sets $E_k$ inductively by $E_0=[0,1]$ and $
E_{k+1}=E_k\setminus \bigcup_c(c+r_{k+1},c+r_k-r_{k+1}), $ where $c$
ranges over the left endpoints of the connected components of $E_k$.
Then $\bigcap_kE_k$ is Cantor. Let $E$ be a $\Q$\nbd linearly
independent Cantor set contained in $\bigcap_kE_k$ (recall the proof
of~\ref{newcantor}). Let $K$ be a Cantor set constructed as
in~\ref{newcantor} with $A=E$; then $(\rbar,K)$ defines a Borel
isomorph of $(\rbar,\N)$. Let $\rr$ be an exponentially bounded
\omin\ expansion of~$\rbar$. We must show that every unary set
definable in $(\rr,K)$ either has interior or is Hausdorff null.
Indeed, with an eye toward applications and further generalization,
we show that every bounded unary set definable in $(\rr,K)^{\#}$
either has interior or is \nul\
(recall~\ref{nullfacts}.\ref{hnull}). We use repeatedly that \nul
ity is preserved under finite unions. By~\ref{FMcase}
and~\ref{newcantor}, it suffices to let $f\colon \R^n\to\R$ be
bounded and definable in $\rr$, and show that $f(E^n)$ is \nul. By
compactness of~$E^n$, it is enough to show that if $x\in E^n$ then
$f(E^n\cap B(x,\delta))$ is \nul\ for some $\delta>0$.
By~\ref{Scellsexp}, there exist $\delta,j>0$ such that for every
$T\in\mathcal T_n$ and $m\leq n$, the restriction of $f$ to $
(x+T(S_{m,j}\times\{0\}^{n-m}))\cap B(x,\delta)$ is continuous and
extends continuously to the closure. Thus, it suffices now
by~\ref{holderimage} to show that:
\begin{itemize}\item[(i)] Every $E^n$ is \nul.
\item[(ii)] For every $j>0$ there
exists $\delta>0$ such that for all $n$,
$$
E^n-E^n\cap (-\delta,\delta)^n\subseteq \bigcup_{m\leq
n}\bigcup_{T\in\mathcal T_n}T(S_{m,j}\times \{0\}^{n-m}),
$$
where $E^n-E^n$ denotes the difference set of $E^n$.
\end{itemize}

\subsubsection*{Proof of~\tx{(i)}} As \nul ity is preserved under finite
cartesian products, it suffices by~\ref{nullfacts}.\ref{ebd} to fix
$m$ and show that $\olim{r}r\exp_m(N(E,r))=0$. Let $0<r<1$; then
there exists $k$ such that $r_{k+1}\leq r<r_k$. For every $j$, the
set $E_j$ consists of $2^{j}$ disjoint closed intervals each of
length $r_j$, so
$$
r\exp_m(N(E,r))\leq r\exp_m(2^{k+1})<r_k\exp_{m+1}(k+1).
%=\frac{r_k}{\psi_{m+1}(\frac{1}{k+1})}.
$$
If $r$ is sufficiently small, then $r_{k-1}\leq 1/(k+1)$ and
$\exp_{m+1}(k+1)=1/\psi_{m+1}(1/(k+1))$, so $r\exp_m(N(E,r))\leq
r_k/\psi_{m+1}(r_{k-1}).$ By construction,
$\ilim{k}r_k/\psi_{m+1}(r_{k-1})=0$.

\subsubsection*{Proof of~\tx{(ii)}} As $j>0$, we have $t-\psi_j(t)>0$ for
all sufficiently small $t>0$. By exponential bounds and properties
of $(r_k)$, we have
$$\ilim{k}\frac{r_{k+1}}{r_k-\psi_j(r_k)}=0=
\ilim{k}\frac{r_{k+1}}{ \psi_j(\psi_j(r_k))}.
$$
Thus, there exists $N\in\N$ such that for all $k>N$ we have
$r_{k+1}<\psi_j(\psi_j(r_k))$ and $\psi_j(r_k)<r_k-r_{k+1}$, hence
also $ r_{k+1}<\psi_j(r_k-r_{k+1}). $ Put $\delta=r_{N}$. We now
proceed by induction on $n$. The case $n\leq 1$ is trivial. Let
$n>1$ and assume the result for all lower values of $n$. The
argument is routine, but a bit tedious to write up in detail; we
give only an outline. By permuting coordinates, it is enough to show
that
$$
E^n-E^n\cap (0,\delta)^n\subseteq \bigcup_{m\leq
n}\bigcup_{T\in\mathcal T_n}T(S_{m,j}\times \{0\}^{n-m}).
$$
By $\Q$\nbd linear independence and symmetry, it is enough to show
that $$ E^n-E^n\cap (0,\delta)^n\cap \mathbb S_n\subseteq S_{n,j}$$
(recall~\ref{phidef} and~\ref{Tdef}).
%Inductively, we have $E^{n-1}-E^{n-1}\cap (0,\delta)^{n-1}\cap
%\mathbb S_{n-1}\subseteq S_{n-1,j}.$
Let $(x,x_{n-1},x_n),(y,y_{n-1},y_n)\in E^{n-2}\times E\times E$ be
such that $(x,x_{n-1},x_n)-(y,y_{n-1},y_n)\in (0,\delta)^n\cap
\mathbb S_n$. Inductively, $(x,x_{n-1})-(y,y_{n-1})\in S_{n-1,j}$,
so it suffices to show that $x_n-y_n<\psi_j(x_{n-1}-y_{n-1})$. Let
$k$ be such that $r_{k+1}<x_{n-1}-y_{n-1}\leq r_k$. It suffices now
by choice of $\delta$ and monotonicity of $\psi_j$ to show that
$x_n-y_n\leq r_{k+1}$. By construction of $E$, we have
$$
E^2\cap [0,r_k]^2\cap \mathbb S_2\subseteq [0,r_{k+1}]^2\cup
[r_k-r_{k+1},r_k]\times [0,r_{k+1}].
$$
%The fiber of $\mathbb S_n$ over $x-x'$ is equal to $\mathbb S_2$.
Since $x_{n-1}-y_{n-1}\in [r_k-r_{k+1},r_k]$, we have $x_n-y_n\in
[0,r_{k+1}]$.
%
%The fiber of $\mathbb S_n$ over points of $\mathbb S_{n-2}$ are
%equal to $\mathbb S_{n-2}$
\end{proof}

Having established Theorem~B, we now proceed to some variants and
corollaries.

Following~\cite{fast}, we say that a sequence $(a_k)$ of positive
real numbers is \textbf{fast} for an expansion $\rr$ of $\rbar$, or
\textbf{$\rr$\nbd fast}, if $\ilim{k}f(a_k)/a_{k+1}=0$ for every
$f\colon\R\to\R$ definable in $\rr$. For $(r_k)$ as in the proof of
Theorem~B, the sequence $(1/r_k)$ is fast for every exponentially
bounded expansion of $\rbar$. An examination of the proof of
Theorem~B yields the following generalization.

\begin{thmc}\label{mainthm}
Let $(a_k)$ be a sequence of positive real numbers. Then there is a
Cantor set $K$ such that $(\R,<,+,1,K)$ $\emptyset$\nbd defines a
Borel isomorph of $(\rbar,\N)$\tx, and for every \omin\ expansion
$\rr$ of $\rbar$\tx, if $(a_k)$ is $\rr$-fast\tx, then every bounded
unary set \dfbl\ in $(\rr,K)^{\#}$ either has interior or is \nul.
\end{thmc}

Under fairly reasonable assumptions, $\rr$-fast sequences exist.

\begin{blank}\label{cntblcof}Let
$\rr$ be an expansion of $\rbar$.
\begin{enumerate}
\item
If there is a countable collection $\mathcal F$ of functions
$\R\to\R$ such that every unary function definable in $\rr$ is
bounded at $+\infty$ by a member of~$\mathcal F$\tx, then there
exist $\rr$\nbd fast sequences.
\item
If $\rr$ is \omin\ and the language of $\rr$ is countable\tx, then
there exist $\rr$-fast sequences.
\end{enumerate}
\end{blank}

\begin{proof}
Diagonalization yields~(1). For~(2), suppose that $\rr$ is \omin. By
the proof of~\cite{geocat}*{C.4}, every unary function definable in
$\rr$ is bounded at $+\infty$ by a unary function \odfbl\ in $\rr$.
Since the language is countable, there are only countably many unary
functions \odfbl\ in $\rr$. Apply~(1).
\end{proof}

When combined with Theorem~C,

\begin{blank}\label{thmagen}
Let $(\rr_k)_{k\in\N}$ be a sequence of \omin\ expansions of
$\rbar$\tx, each in a countable language. Then there is a Cantor set
$K$ such that $(\rbar,K)$ defines a Borel isomorph of $(\rbar,\N)$
and every bounded unary set definable in any $(\rr_k,K)^{\#}$ either
has interior or is Minkowski null.
\end{blank}

\begin{blank}\label{discretecor}
Theorem~C holds with \lqu discrete\rqu\tx, hence also \lqu
countable\rqu\tx, in place of \lqu Cantor\rqu.
\end{blank}

\begin{proof}
Let $E$ be any Cantor set and $M$ be the set of midpoints of the
(bounded) complementary intervals of~$E$. Note that
$$M=\set{r\in \R: \exists\e>0, E\cap
[r-\e,r+\e]=\{r-\e,r+\e\}}$$ and $E=\Fr(M)$. Hence, $(\rbar,E)$ and
$(\rbar,M)$ are $\emptyset$\nbd interdefinable.
\end{proof}

We now answer a question raised in~\cite{tameness}*{\S3.1}. A set
$A\subseteq \R^n$ \textbf{has a locally closed point} if it has
nonempty interior in its closure. It is easy to see that if every
nonempty unary set definable in $\rr$ has a locally closed point,
then every definable unary set either has interior or is nowhere
dense. We show that the converse fails.

\begin{blank}\label{nolocc}
There exist $\emptyset\neq A\subseteq\R$ having no locally closed
points such that every unary set \dfbl\ in $(\rbar,A)$ either has
interior or is nowhere dense.
\end{blank}

\begin{proof}
With $K$ as in Theorem~B, let $A$ be the set of left endpoints of
the complementary intervals of $K$. Note that $A\subseteq K$ and
$\Cl(A)=\Cl(K\setminus A)=K$.
\end{proof}

\section{Discussion and open issues}\label{S:conc}

%\begin{dblank}
It is easy to construct Cantor sets $E$ that do define $\N$ over
$\rbar$: Just encode $\N$ by the set of lengths of the complementary
intervals. For example, remove successively from $[0,1]$ the middle
intervals of length $1/(n!)$ for $n\geq 2$. Then $(\rbar,E)$ defines
the set $A:=\set{n!: n\in \N}$, hence also the successor function
$\sigma\colon A\to A$, hence also $\N=\{0\}\cup \set{\sigma(a)/a:
a\in A}$.
%\end{dblank}

%\begin{dblank}
Though Theorem~A answers one question, others arise immediately. Are
there closed $E\subseteq\R$ and $N\in\N$ such that $(\rbar,E)$
defines a non-Borel set, yet every set definable in $(\rbar,E)$ is
$\bsig_N^1$? (If so, then $N\geq 2$ by Souslin's
Theorem~\cite{kechris}*{14.11}.) Evidently, one can generalize this
question. For example, allow $E$ to be constructible, or $\fsig$, or
of finite Borel rank. We can go in the other direction as well. In
Theorem~A, can we take $E$ to be closed and countable? Closed and
discrete? (\cf~\ref{discretecor}.)
%Contained in $\N$?
Regarding the last, there are some known restrictions. If
$E\subseteq (0,\infty)$ is infinite, closed and discrete, then it is
the range of a strictly monotone and unbounded-above
sequence~$(a_k)$ of positive real numbers. If $(\log a_{k+1})/(\log
a_k)\to \infty$, that is, if $(a_k)$ is $\rbar$\nbd fast, then every
set definable in $(\rbar,E)$ is constructible~\cites{fast}. On the
other hand, if $(a_k)=(f(k))$ for some sufficiently well behaved
$f\colon \R\to\R$---in particular, if $(\rbar,f)$ is \omin---and
$(\log a_k)/k\to 0$, then $(\rbar,E)$ defines $\N$. See~\cite{aph}
for a proof of the latter and some information on behavior between
these extremes.
%A.~Dolich has pointed out to us in correspondence that if
%$E\subseteq \N$ and $(\rbar,E)$ does not define $\N$, then $E$ has
%asymptotic density zero. To see this, suppose otherwise. The set
%$$
%A:=\set{x\in\R: \exists a,b\in E,\ b\neq 0\And \ a+xb\in E}
%$$ is definable in $(\rbar,E)$ and contained in $\Q$.
%By Szemer\'edi's Theorem~\cite{szem}, $E$ contains arbitrarily long
%arithmetic sequences, so $\N\subseteq A$ and $\Q=\set{\pm x/y:
%x,y\in A\And y\neq 0}$. Thus, $(\rbar,E)$ defines $\Q$, hence also
%$\N$ by J.~Robinson~\cite{jrob}.

(There are non-Borel $E\subseteq \R$ such that every set definable
in $(\rbar,E)$ is $\bsig_2^1$, but this is far off the point of this
paper, so we give only a hint: By van den
Dries~\cite{densepairs}*{Theorem~1}, if $E$ is $\bsig_N^1$ and a
real-closed subfield of~$\R$, then every set definable in
$(\rbar,E)$ is a boolean combination of $\bsig_{N}^1$ sets.)
%\end{dblank}

%\begin{dblank}
Currently, we know of no expansions of $\rbar$ by constructible
subsets of $\R$ that define non-constructible sets, yet every
definable set is Borel. There are candidates, though, as we explain
in the next two paragraphs.

The proof of~\ref{nolocc} shows that if $E$ is a Cantor set, then
$(\rbar,E)$ defines a unary set that is not $\fsig$ (hence not
constructible). Are there Cantor sets $E$ such that every set
definable in $(\rbar,E)$ is a boolean combination of $\fsig$ sets?
Along these lines, the Cantor set $K$ in~\ref{newcantor} was
designed to encode $(\rbar,\N)$. But we could omit deliberately
encoding $\N$, that is, in the proof of~\ref{newcantor}, replace $X$
with its projection on the first four coordinates. Define $E$ as in
the proof of Theorem~B, and so on. What can be said about the
definable sets of $(\rbar,K)$? (Of course, by the rest of the proof
of Theorem~B, every unary definable set either has interior or is
Hausdorff null.) Or, in the definition of $X$, replace $\dis_{\N}$
with $\exp$. Then $(\rbar,K)$ Borel-interprets $(\rbar,\exp)$. Is
$\N$ definable? Is $\exp$? And so on. The technology
from~\cite{ominsparse} appears not to tell us much more about the
sets definable in $(\rbar,E)$ than those in $(\rbar,E)^{\#}$, so it
seems that new ideas are needed.

Given $\alpha>0$, it is known that every set definable in
$(\rbar,\alpha^{\N})$ is constructible~\cite{opencore}*{\S4} (and
more~\cites{ominsparse,tameness}), where
$\alpha^{\N}=\set{\alpha^n:n\in\N}$. The expansion of $\rbar$ by
the set $\set{2^n+3^n:n\in\N}$ defines both $2^\N$ and $3^\N$
(see~\cite{aph}), hence also the multiplicative group $2^{\Z}\cdot
3^{\Z}$, which is dense and co-dense in $(0,\infty)$, hence not
constructible. Aside from defining $2^{\Z}\cdot 3^{\Z}$, little is
known about $(\rbar,2^{\N},3^{\N})$, in particular, whether it
defines any sets that are not ${\boldsymbol\Delta}_3^0$. There is
nothing special here about the bases $2$ and $3$: the situation is
the same for any $\alpha_1,\dots,\alpha_N>0$ such that $\log
\alpha_1\dots,\log\alpha_N$ are $\Q$\nbd linearly independent.

We close with a few words on history and attribution. Recall that
our original goal was to find a constructible $E\subseteq \R$ such
that $(\rbar,E)$ defines a non-Borel set but does not define $\N$.
Every expansion of $\rbar$ that defines $\N$ also $\emptyset$\nbd
defines it, hence also~$\Q$, because $\N$ is the unique subset of
$[0,\infty)$ that is closed under $x\mapsto x+1$ and whose
intersection with $[0,1]$ is equal to~$\{0,1\}$. Hence, the
following \textit{a priori} weaker version of Theorem~A due to
Friedman and Miller suffices.
\begin{thma0}
There is a closed $E\subseteq\R$ such that $(\rbar,E)$ defines a
Borel isomorph of $(\rbar,\N)$ and every unary set \odfbl\ in
$(\rbar,E)$ either has interior or is nowhere dense.
\end{thma0}
\noindent For this, it suffices by~\ref{sparsefact}
and~\ref{newcantor} to find a Cantor set $E$ such that $f(E^n)$ is
nowhere dense for every $f\colon\R^n\to\R$ that is \odfbl\ in
$\rbar$. The naive approach described in the introduction appears to
stall for functions \odfbl\ in $\rbar$ just as it does for
parametrically definable functions. But in any \omin\ structure in a
countable language, there are only countably many \odfbl\ functions,
and they are all Borel (by cell decomposition). This prompted
Friedman to produce the following result of independent interest,
thus establishing an appropriately modified version
of~\ref{thmagen}, hence also Theorem~A${_0}$.

\begin{blank}\label{borelfcns}
For every sequence $(f_k\colon \R^{n(k)}\to\R)_{k\in\N}$ of Borel
functions there is a Cantor set $E$ such that every image
$f_k(E^{n(k)})$ is nowhere dense.
\end{blank}
\noindent We shall not prove this here as we no longer need it; the
interested reader may wish to attempt verification by amalgamating
the proofs of~19.1 and~19.8 from~\cite{kechris}. On the other hand,
attempts by Friedman and Miller to derive Theorem~A from
Theorem~A${_0}$ were unsuccessful, as were attempts to conclude
from~\ref{borelfcns} the existence of a Cantor set $E$ such that
$f(E^n)$ is nowhere dense for every semialgebraic
$f\colon\R^n\to\R$. Of course, \ref{borelfcns} is a very blunt
hammer in this setting, as it uses nothing about \omin ity (yet
relies heavily on countability). Convinced that a more
singularity-theoretic approach was in order, Miller approached
Kurdyka and Speissegger, who subsequently solved the semialgebraic
case, which Miller then refined to its current form. The crucial
idea of using Minkowksi rather than Hausdorff nullity is due to
Kurdyka. Result~\ref{bddcell} is due to Speissegger, who had known
it for several years but not made prior use of it.
%\end{dblank}

\bibsection

\end{document}